\documentclass{amsart}

\usepackage[utf8]{inputenc}
\usepackage[T1]{fontenc}
\usepackage{lmodern}
\usepackage[french]{babel}
\usepackage{amsmath}
\usepackage{amsthm}
\usepackage{amssymb}
\usepackage{units}
\usepackage{amscd}
\usepackage{geometry}
\usepackage{graphicx}
\usepackage{tikz}
\usepackage{tikz-cd}
\usepackage{pst-plot}
\usepackage{enumitem}
\usetikzlibrary{arrows,automata}
\usepackage{graphicx}
\usepackage{url}

{\newtheorem{theo}{Théorème}[section]}
{}
{}
{\theoremstyle{definition} \newtheorem{defin}[theo]{Définition}
							\newtheorem{prop}[theo]{Proposition}
							\newtheorem{lemme}[theo]{Lemme}}
{\theoremstyle{remark} \newtheorem{remarque}[theo]{Remarque}
						}

	\newcommand{\Id}{Id}	
\newcommand{\Set}{\text{Set}}
\newcommand{\op}{op}

\newcommand{\sS}{\text{sSet}}

\newcommand{\Hom}{\text{Hom}}

\newcommand{\Top}{\text{Top}}
\newcommand{\Real}[1]{{}||#1||}
\newcommand{\RealP}[1]{{}\Real{#1}_P}
\newcommand{\Sing}{\text{Sing}}

\newcommand{\pr}{pr}

\newcommand{\Fun}{\text{Fun}}

\newcommand{\ev}{\text{ev}}

\newcommand{\Hol}{\text{Holink}}

\newcommand{\fil}[1]{(#1,\varphi_{#1})}

\newcommand{\Strat}{\text{Strat}}

\title{Théorème de Whitehead stratifié et invariants de noeuds}
\author{Sylvain Douteau}

\begin{document}
\begin{abstract}
En considérant des homotopies préservant la stratification, on obtient une notion naturelle d'homotopie pour les espaces stratifiés. Dans cette note, on présente des invariants d'homotopie stratifiée, les groupes d'homotopie stratifiés. On montre que ces groupes d'homotopie stratifiés vérifient un analogue stratifié au théorème de Whitehead. Comme illustration, on présente un invariant de nœud complet défini à partir des groupes d'homotopies stratifiés.
\end{abstract}
\maketitle

\section{Espaces stratifiés}
\begin{defin}
Un espace stratifié est une paire $\fil{X}$ où
\begin{itemize}
\item $X$ est un espace topologique,
\item $\varphi_X\colon X\to P$ est une application continue, entre $X$ et un ensemble ordonné $P$ muni de la topologie d'Alexandroff.
\end{itemize}
Une application stratifiée $f\colon\fil{X}\to\fil{Y}$ est la donnée d'une paire d'applications continues $f,\widehat{f}$ telles que le diagramme suivant commute
\begin{equation*}
\begin{tikzcd}
X
\arrow{r}{f}
\arrow[swap]{d}{\varphi_X}
&Y
\arrow{d}{\varphi_Y}
\\
P
\arrow[swap]{r}{\widehat{f}}
&Q
\end{tikzcd}
\end{equation*}
Ceci définit la catégorie des espaces stratifiés $\Strat$.
\end{defin}

\begin{remarque}\label{PseudoVariete}
Si $X$ est une pseudo-variété (voir par exemple \cite[Definition 2.4.1]{Friedman}), on considère $P_X$ l'ensemble des strates de $X$. Alors, on a une relation d'ordre sur $P_X$ donnée par $S\leq T\Leftrightarrow S\subset \overline{T}$. On vérifie que l'application $\varphi_X\colon X\to P_X$, qui à tout point de $X$ associe la strate à laquelle il appartient, est continue. En particulier, la catégorie des pseudo-variétés et des morphismes préservant les strates forme une sous catégorie pleine de $\Strat$.
\end{remarque}

\begin{defin}
Soient $f,g\colon \fil{X}\to\fil{Y}$ deux applications stratifiées. On dit que $f$ et $g$ sont homotopes par une homotopie stratifiée si il existe une application stratifiée $H$ 
\begin{equation*}
\begin{tikzcd}
X\times [0,1]
\arrow[swap]{d}{\varphi_X\circ\pr_X}
\arrow{r}{H}
&Y
\arrow{d}{\varphi_Y}
\\
P
\arrow[swap]{r}{\widehat{H}}
&Q
\end{tikzcd}
\end{equation*}
telle que $H_{|X\times\{0\}}=f$ et $H_{|X\times\{1\}}=g$.
On dit que $f\colon \fil{X}\to\fil{Y}$ est une équivalence d'homotopie stratifiée si il existe une application stratifiée $g\colon\fil{Y}\to\fil{X}$ telle qu'il existe des homotopies stratifiées entre $f\circ g$ et $\Id_Y$ et entre $g\circ f$ et $\Id_X$.
\end{defin}

\begin{remarque}\label{PosetConstant}
Si $f$ et $g$ sont deux applications stratifiées homotopes par une homotopie stratifiée, alors $\widehat{f}=\widehat{g}$. En particulier, si $f$ est une équivalence d'homotopie stratifiée, $\widehat{f}$ est un isomorphisme d'ensemble ordonnés. Pour la suite de ce texte, on fixe un ensemble ordonné $P$, et on considère des applications stratifiées $f$ telle que $\widehat{f}=\Id_P$. On note $\Top/P$ la sous catégorie de $\Strat$ correspondante. 
\end{remarque}

\begin{defin}
On définit la catégorie $\Delta(P)$ des $P$-simplexes filtrés comme suit (voir \cite{MonArticle1}):
\begin{itemize}
\item Les objets sont les morphismes d'ensembles simpliciaux $\varphi\colon\Delta^n\to N(P)$, où $\Delta^n$ est un simplexe et $N(P)$ est le nerf de l'ensemble ordonné $P$.
\item Les morphismes sont donnés par 
\begin{equation*}
\Hom(\varphi\colon\Delta^n\to N(P),\psi\colon\Delta^m\to N(P))=\{f\colon \Delta^n\to \Delta^m\ |\ \psi\circ f=\varphi\}
\end{equation*}
\end{itemize}
On définit la catégorie $R(P)$ des $P$-simplexes réduits comme la sous catégorie pleine de $\Delta(P)$ contenant les objets tels que $\varphi\colon\Delta^n\to N(P)$ est un monomorphisme. Pour simplifier les notations, on note $\Delta^{\varphi}$ le simplexe filtré $\varphi\colon\Delta^n\to N(P)$.
\end{defin}

\begin{remarque}
Par des arguments généraux, on a l'équivalence de catégorie
\begin{equation*}
\Fun(\Delta(P)^{\op},\Set)\simeq\sS/N(P)
\end{equation*}
où $\sS/N(P)$ est la catégorie des ensembles simpliciaux au dessus de $N(P)$.
\end{remarque}

\begin{prop}\label{Adjonction}
Il existe une adjonction
\begin{equation*}
\RealP{-}\colon\sS/N(P)\leftrightarrow \Top/P\colon \Sing_P
\end{equation*}
où $\RealP{K\to N(P)}$ est donné par la composition
\begin{equation*}
\Real{K}\to\Real{N(P)}\to P
\end{equation*}
Voir \cite[Definition A.6.2 et Remark A.6.3]{HigherAlgebra}
\end{prop}

\begin{defin}
Soit $\fil{X}$ un espace stratifié. Un pointage de $\fil{X}$ est la donnée d'un sous ensemble simplicial $V\subseteq N(P)$ et d'une application stratifiée 
\begin{equation*}
\phi\colon\RealP{V}\to\fil{X}
\end{equation*}
\end{defin}

\begin{defin}
Soit $\fil{X}$ un espace stratifié, et $\phi$ un pointage de $X$. On définit le foncteur d'entrelacs homotopiques de $\fil{X}$ comme
\begin{align*}
\Hol(\fil{X},\phi)\colon R(P)^{\op}&\to\Top_*\\
\Delta^{\varphi}&\mapsto \left\{\begin{array}{cl}
\left(C^0_P\left(\RealP{\Delta^{\varphi}},\fil{X}\right)\phi_{|\Delta^{\varphi}}\right)&\text{ Si $\Delta^{\varphi}\subseteq V$}\\
(*,*)&\text{ sinon}
\end{array}\right. 
\end{align*}
Ici, $C^0_P\left(\RealP{\Delta^{\varphi}},\fil{X}\right)$ est l'ensemble $\Hom_{\Top_P}(\RealP{\Delta^{\varphi}},\fil{X})$ muni de la topologie induite par l'inclusion
\begin{equation*}
\Hom_{\Top_P}(\RealP{\Delta^{\varphi}},\fil{X})\subseteq C^0(\Real{\Delta^{n}},X).
\end{equation*}
En composant avec le foncteur $\pi_n$, on obtient les groupes d'homotopies stratifiés de $\fil{X}$
\begin{equation*}
s\pi_n(\fil{X},\phi)=\pi_n(\Hol(\fil{X},\phi))
\end{equation*}
\end{defin}

\section{Théorème de Whitehead stratifié}
\begin{theo}\label{WhiteheadFiltre}
Soient $X$ et $Y$ deux pseudo-variétés PL, et $f\colon X\to Y$ une application continue préservant les strates. Alors, $f$ admet un inverse à homotopie stratifiée près si et seulement si, pour tout pointage $\phi\colon\RealP{V}\to\fil{X}$ et pour tout $n\geq 0$, $f$ induit des isomorphismes
\begin{equation*}
s\pi_n(\fil{X},\phi)\to s\pi_n(\fil{Y},f\circ\phi)
\end{equation*}
\end{theo}

Par la remarque \ref{PseudoVariete}, $X$ et $Y$ sont des espaces stratifiés au dessus de $P_X$ et $P_Y$, et par la remarque \ref{PosetConstant} on se ramène au cas où $P_X=P_Y$. 
La preuve de ce théorème découle des résultats suivants :

\begin{theo}[{\cite[Theorem 3.19]{MonArticle1}}]
Il existe une structure de catégorie modèle sur $\sS/N(P)$ où les cofibrations sont les monomorphismes et où les équivalences d'homotopie sont les équivalences d'homotopie stratifiées. De plus, les équivalences faibles entre objets fibrants sont exactement les morphismes induisant des isomorphismes entre groupes d'homotopies stratifiés.
\end{theo}

\begin{lemme}[{\cite[Theorem A.6.4]{HigherAlgebra}}]
Si $X\to P_X$ est une pseudo-variété stratifiée comme dans la remarque \ref{PseudoVariete}, $\Sing_P(X)$ est fibrant.
\end{lemme}

\begin{lemme}
L'adjonction de la proposition \ref{Adjonction} préserve les homotopies stratifiés. De plus, si $X$ est une pseudo-variété, pour tout pointage $\phi$ de $X$, $s\pi_n(X,\phi)=s\pi_n(\Sing_P(X),\phi)$.
\end{lemme}

\begin{proof}[Démonstration du théorème \ref{WhiteheadFiltre}]
Le sens direct du théorème \ref{WhiteheadFiltre} est une conséquence de la définition de $s\pi_n$.
L'idée de la preuve de la réciproque est contenue dans le diagramme commutatif suivant. 
\begin{equation}\label{DiagrammePreuveWhitehead1}
\begin{tikzcd}[column sep = huge]
\phantom{X}
&\RealP{\Sing_P(Y)}
\arrow{r}{\RealP{\widetilde{g}}}
&\RealP{\Sing_P{X}}
\arrow{r}{\RealP{\Sing_P(f)}}
\arrow{d}{\ev_X}
&\RealP{\Sing_P{Y}}
\arrow{d}{\ev_Y}
\\
X
\arrow{d}{i_X}
\arrow{r}{h}
&Y
\arrow{r}{g}
\arrow{u}{i_Y}
&X
\arrow{r}{f}
&Y
\\
\RealP{\Sing_P{X}}
\arrow{r}{\RealP{\widetilde{h}}}
&\RealP{\Sing_P{Y}}
\arrow{u}{\ev_Y}
\arrow{r}{\RealP{\Sing_P(g)}}
&\RealP{\Sing_P{X}}
\arrow{u}{\ev_X}
\end{tikzcd}
\end{equation}
On suppose que $f$ induit des isomorphismes entre les groupes d'homotopie stratifiés. C'est donc aussi le cas pour $\Sing_P(f)$. On en déduit que $\Sing_P(f)$ est une équivalence faible entre objets fibrants et cofibrants donc une équivalence d'homotopie stratifiée. On note $\widetilde{g}$ un inverse à homotopie stratifiée près. Comme $Y$ est PL, on dispose d'une application stratifiée canonique $Y\to \RealP{\Sing_P(Y)}$, ce qui permet de définir $g$. Comme $\widetilde{g}$ est un inverse de $\Sing_P(f)$, et que l'adjonction $\RealP{-},\Sing_P$ préserve les homotopies stratifiées, on a que $g$ est un inverse à gauche de $f$ à homotopie stratifiée près. On en déduit en particulier que $g$ induit des isomorphismes sur tout les groupes d'homotopie stratifiés. En répétant la construction pour $g$, on obtient $h$ un inverse à gauche de $g$ à homotopie stratifiée près. Puisque $g$ admet un inverse à gauche et à droite, c'est une équivalence d'homotopie stratifiée, d'inverse $f$. 
\end{proof}

\section{Exemple}

\begin{defin}
Soit $\gamma\colon S^1\to S^3$ un noeud lisse. On note $K=\gamma(S^1)$. On voit $S^3$ comme un espace stratifié au dessus de $P=\{0<1\}$ en définissant $\varphi\colon S^3\to P$ comme 
\begin{equation*}
\varphi(x)=\left\{\begin{array}{cl}
0 & \text{ si $x\in K$}\\
1 & \text{ si $x\in S^3\setminus K$}
\end{array}\right.
\end{equation*} 
On note $X_{\gamma}$ l'espace stratifié correspondant.
\end{defin}

\begin{remarque}
L'espace stratifié ainsi obtenu correspond à une pseudo-variété.
\end{remarque}

\begin{remarque}
l'espace topologique $C^0_P([0,1],X_{\gamma})$ est l'espace des applications $f\colon [0,1]\to S^3$ telles que $f(0)\in K$, et $f(t)\in S^3\setminus K$, $t>0$. On vérifie facilement qu'il est faiblement équivalent au bord d'un voisinage tubulaire de $K$, c'est à dire à un tore.
\end{remarque}

Fixons un pointage de $X_{\gamma}$, $\phi\colon [0,1]\to X_{\gamma}$, et calculons son premier groupe d'homotopie stratifié. D'après la remarque précédente, on obtient le diagramme suivant :

\begin{equation*}
s\pi_1(X_{\gamma},\phi)=
\begin{tikzcd}
&\mathbb{Z}\oplus\mathbb{Z}
\arrow[swap]{dl}{\pr_1}
\arrow{dr}{f}
\\
\mathbb{Z}
&
&
G_{\gamma}
\end{tikzcd}
\end{equation*}
où $G_{\gamma}$ est le groupe fondamental du complémentaire du nœud, et $f$ est induit par l'inclusion du bord d'un voisinage tubulaire de $K$ dans le complémentaire de $K$. En particulier, $f(\mathbb{Z}\oplus\mathbb{Z})\subset G_{\gamma}$ est le sous-groupe périphérique. Finalement, on a le résultat suivant

\begin{theo}
Soient $\gamma_1,\gamma_2\colon S^1\to S^3$ deux noeuds lisses. Il existe un homéomorphisme $f\colon S^3\to S^3$ tel que $f(\gamma_1(S^1))=\gamma_2(S^1)$ si et seulement si il existe une équivalence d'homotopie stratifiée $g\colon X_{\gamma_1}\to X_{\gamma_2}$.
\end{theo}

\begin{proof}
Tout homéomorphisme vérifiant les hypothèses du théorème précédent est en particulier une équivalence d'homotopie stratifiée. Réciproquement, une équivalence d'homotopie stratifiée induit un isomorphisme entre les groupes de noeuds préservant les sous-groupes périphériques. Par les résultats de \cite{Waldhausen} et de \cite{GordonLuecke} cela implique qu'il existe un homéomorphisme vérifiant les hypothèses du théorème.
\end{proof}

\bibliographystyle{alpha}
\bibliography{biblio}

\begin{thebibliography}{{Dou}18}

\bibitem[{Dou}18]{MonArticle1}
Sylvain {Douteau}.
\newblock {A Simplicial Approach to Stratified Homotopy Theory}.
\newblock {\em ArXiv e-prints}, January 2018.

\bibitem[Fri]{Friedman}
Greg Friedman.
\newblock {\em Singular Intersection Homology}.
\newblock \url{http://faculty.tcu.edu/gfriedman/IHbook.pdf}.

\bibitem[GL89]{GordonLuecke}
C.~McA. Gordon and J.~Luecke.
\newblock Knots are determined by their complements.
\newblock {\em J. Amer. Math. Soc.}, 2(2):371--415, 1989.

\bibitem[Lur]{HigherAlgebra}
Jacob Lurie.
\newblock {\em Higher Algebra}.
\newblock \url{http://www.math.harvard.edu/~lurie/papers/HA.pdf}.

\bibitem[Wal68]{Waldhausen}
Friedhelm Waldhausen.
\newblock On irreducible {$3$}-manifolds which are sufficiently large.
\newblock {\em Ann. of Math. (2)}, 87:56--88, 1968.

\end{thebibliography}

\end{document}